\documentclass[11pt]{article}
\usepackage{amssymb}
\usepackage{}
\usepackage{stmaryrd}
\usepackage{amsmath,amsfonts}
\usepackage{color,txfonts,wasysym,lscape,amssymb,mathrsfs,extarrows}
\usepackage[pagewise]{lineno}
\usepackage[all,pdf]{xy}

\newtheorem{theorem}{Theorem}[section]
\newtheorem{proposition}{Proposition}[section]

\newtheorem{definition}{Definition}[section]

\newtheorem{remark}{Remark}[section]

\title{Q-Linear Convergence of the Proximal Augmented Lagrangian Method for Non-Convex Conic Programming}

 \author{Ning Zhang \footnote{ School of Computer Science and Technology, Dongguan University of Technology, Dongguan, 523808, China.
 		(e-mail:zhangning@dgut.edu.cn). This work was supported by the National Natural Science Foundation of China under grant number 12271095 and 121701185.}\quad Yi Zhang \footnote{Corresponding author. School of Mathematics, East China University of Science and Technology,
 Shanghai, 200237, China.(e-mail:yizhang@ecust.edu.cn).  This work was supported by the National Natural Science Foundation of China under grant number 12171153.}
}

\date{\today}
\begin{document}

\maketitle

\begin{abstract} 	
This paper provides a local convergence analysis of the proximal augmented Lagrangian method (PALM) applied to a class of non-convex conic programming problems. Previous convergence results for PALM typically imposed assumptions such as constraint non-degeneracy, strict complementarity, second-order sufficiency conditions, or a combination of constraint non-degeneracy with strong second-order sufficiency conditions. In contrast, our work demonstrates a Q-linear convergence rate for an inexact version of PALM in the context of non-convex conic programming, without requiring the uniqueness of the Lagrange multipliers. The analysis relies solely on the second-order sufficiency condition and the calmness property of the multiplier mapping, presenting a more relaxed set of conditions for ensuring convergence.

\vskip 12 true pt \noindent \textbf{Key words}: the proximal augmented Lagrangian method, calmness property, second-order sufficient condition, conic programming, Q-linear convergence.

\vskip 12 true pt \noindent \textbf{MSC 2010 }: 9C30
\end{abstract}

\section{\textbf{Introduction}}
 \setcounter{equation}{0}

Let $X$ and $H$ be two finite-dimensional Euclidean spaces, and consider the constrained optimization problem:
\begin{align}\label{1.1}
	\min_{x\in X} \ f(x)\quad \mbox{s.t.} \ g(x)\in K,
\end{align}
where $f:X\rightarrow \mathbb{R}$ and $g:X\rightarrow H$ are both twice  differentiable functions, and $K\subseteq H$ is a non-empty closed convex set. In this paper, we assume that $K$ is $\mathcal{C}^2-$cone reducible (as defined in Definition \ref{d2.2}). Under this assumption, problem (\ref{1.1}) includes nonlinear programming, semidefinite programming, second-order cone programming, and any combination of these \cite{Shapiro03}.

 The augmented Lagrangian $\mathcal{L}:X\times H\times (0,+\infty)\rightarrow \mathbb{R}$ associated with (\ref{1.1}) is defined by:
\begin{align}\label{augmented function}
	\mathcal{L}(x,\lambda,c):=f(x)+\frac{c}{2}\mbox{dist}^2\left(g(x)+c^{-1}\lambda,K\right)-\frac{1}{2}c^{-1}\|\lambda\|^2,
\end{align}
where $\lambda$ is a multiplier and $c>0$ is a penalty parameter. Given the current iterate $(x_k, \lambda_k, c_k, \epsilon_k)$, where $x_k\in X$ is the primal variable, $\lambda_k \in H$ is the dual variable, $c_k > 0$ is the penalty parameter and $\epsilon_k \geq 0$ is the error tolerance parameter.
The classical augmented Lagrangian method (ALM), which is inexact when  $\epsilon_k>0$, generates the next primal-dual iterate as follows.
The primal iterate \(x^{k+1} \in X\) is a approximate solution of the unconstrained optimization problem:
\begin{align}\label{ALM3}
	\min_{x \in X} \mathcal{L}(x, \lambda^k, c_k),
\end{align}
where \(x^{k+1}\) satisfies the first-order optimality condition:
\begin{align}\label{ALM1}
	\|\nabla_x \mathcal{L}(x^{k+1}, \lambda^k, c_k)\| \leq \epsilon_k.
\end{align}
The update of the dual iterate is given by
\begin{align}\label{lambda-update}
	\lambda^{k+1}=c_kg(x^{k+1})+\lambda^k-\Pi_{K}\left(c_kg(x^{k+1})+\lambda^k\right).
\end{align}

The literature on the ALM is extensive, and we focus primarily on studies related to local convergence results and convergence rates. This method was initially proposed by Hestenes \cite{hestenes1969} and Powell \cite{powell1969} for equality-constrained problems, and later extended to nonlinear programming (NLP) by Rockafellar \cite{Rockafellar1973}. A considerable body of theoretical work has been dedicated to the convergence analysis of ALM, including foundational studies by \cite{bertsekas1982, golshtein1989, Rockfellar1976, rockafellar1993}, among others.
For non-convex NLP problems, Bertsekas \cite{bertsekas1982} proved that under the strong second-order sufficient condition (SOSC), linear independence constraint qualification (LICQ), and strict complementarity conditions, the sequence of dual variables converges Q-linearly, while the sequence of primal variables converges R-linearly. Subsequent research focused on relaxing these conditions: Conn et al. \cite{conn1991}, Contesse-Becker \cite{contessebecker1993}, and Ito and Kunisch \cite{ito1990} established linear convergence rates of ALM for general NLP problems without the strict complementarity condition. Further advancements were made by Fernandez and Solodov \cite{fernandez2012}, who removed the uniqueness requirement for Lagrange multipliers, marking a significant breakthrough in convergence theory.
Recent advancements include: Hang and Sarabi \cite{Hang2021}, who proved the local convergence of piecewise quadratic composite optimization problems under the SOSC condition alone; Sun et al. \cite{sun2008}, who established the convergence rate analysis for nonlinear semidefinite programming (NLSDP) under the strong SOSC and constraint non-degeneracy conditions. Kanzow and Steck \cite{kanzow2019} proved the primal-dual linear convergence properties of ALM for constrained optimization problems that satisfy SOSC and strong Robinson constraint qualifications (including NLSDP and nonlinear second-order cone programming (NLSOCP)). Further research on the linear convergence of NLSOCP and NLSDP, as shown in the works of \cite{hang2021} and \cite{Wang2024}, primarily relies on SOSC and the semi-isolated calmness property of the KKT solution map. In recent advancements, Hang \cite{Hang2023} established local convergence theory for the ALM applied to composite optimization problems with non-unique Lagrange multipliers, under the SOSC. This work extends the traditional convergence analysis framework and provides theoretical guarantees for a broader class of optimization problems.

All the research findings mentioned above focus on the basic ALM (\ref{ALM1}) and  (\ref{lambda-update}). The objective of this paper is to extend these theoritical frameworks to an important variant:
the proximal augmented Lagrangian method (PALM), which is proposed in  \cite{Rockfellar1976} for convex programming. For the PALM, (\ref{ALM3}) is replaced by:
\begin{align}\label{x-update}
	\min_{x\in X} \mathcal{L}(x,\lambda^k,c_k) + \frac{1}{2c_k} \|x - x^k\|^2,
\end{align}
and accordingly, (\ref{ALM1}) is replaced by:
\begin{align}\label{x-pupdate}
	\|\nabla_x \mathcal{L}(x^{k+1}, \lambda^k, c_k) + c_k^{-1} (x^{k+1} - x^k)\| \leq \epsilon_k,
\end{align}
where $c^{k+1}$ is updated by a certain rule. PALM extends ALM by incorporating a proximal term, which enhances stability in the optimization process.
 The proximal term improves the method's ability to deal with non-smooth optimization problems, ensuring better convergence and more reliable performance. PALM also adapts more effectively to large-scale problems and complex constraints, making it a more robust choice for challenging optimization tasks.

It is important to note that extending the linear convergence results of the classical ALM to its proximal variants is neither straightforward nor automatically guaranteed. This is because the primal step in (\ref{x-update}) is shorter than that in (\ref{ALM3}). The PALM for convex problems was first studied by Rockafellar \cite{Rockfellar1976}, which can also be be interpreted  as a primal-dual inexact proximal point method. However, this method has received limited attention for non-convex problems until recently.
Zhang et al. \cite{zhang2019} studied a PALM for equality-constrained optimization problems and achieved linear or superlinear convergence rates. PALM has also been applied to NLSDP and NLSOCP in \cite{Chu2021, zhang2020}. The linear and superlinear rates in \cite{Chu2021, zhang2019,zhang2020} were derived based on the same strong assumptions typically used for the ALM, such as constraint non-degeneracy conditiosn (LICQ for NLP) , strict complementarity, and the SOSC (or a combination of constraint non-degeneracy conditiosn and strong SOSC). More recently, Izmailov and Solodov \cite{IS2024} introduced a perturbed ALM framework and established the convergence properties of PALM for NLPs under weaker conditions: once differentiability of the functions and SOSC (without requiring LICQ or strict complementarity). This motivates us to investigate whether the assumptions in the works of \cite{Chu2021, zhang2019, zhang2020} for non-polyhedral problem can be relaxed.  In particular, we aim to establish a locally Q-linear convergence of the PALM for general non-convex $\mathcal{C}^2$-cone reducible problems only under SOSC and the calmness property of the multipler mapping, without assuming the uniqueness of Lagragngian multipliers.

The remaining part of this paper is structured as follows. Section 2 presents the necessary preliminaries, including optimality conditions, local error bounds relevant to conic programming, and an iterative framework for solving generalized equations. Section 3 verifies the solvability of the subproblem in (\ref{x-update}) and presents the key result regarding the local Q-linear convergence of PALM for conic programming, based on the iterative framework for generalized equations. Finally, Section 4 concludes the paper.

\section{Notations and preliminaries }
\setcounter{equation}{0}
We begin by introducing the notations that will be used throughout the paper. Let
 \( X \) represent a finite-dimensional Euclidean space.  The inner product in  \( X \) by \( \langle \cdot, \cdot \rangle \), and the corresponding induced norm is denoted by \( \|\cdot\| \).
 The unit ball centered at the origin in \( X \) is denoted by \( \mathbb{B}_X \). For a given point \( x \in X \), \( B_{\varepsilon}(x) \) denotes the closed ball of radius \( \varepsilon \) centered at \( x \). For a closed convex set \( S \) in \( X \) and a point \( x \) belong to \( X \), the distance from \( x \) to \( S \) is given by $\text{dist}(x, S) := \inf_{y \in S} \|y - x\|$. The projection of \( x \) onto \( S \) is denoted by \( \Pi_S(x) \). For a differentiable function $g$, we denote the derivative of $g$ at $x$ by $\mathcal{J}g(x)$ and denoted the adjoint of $\mathcal{J}g(x)$ by $\nabla g(x)$. We denoted by $\mathbb{R}_+$ the set of non-negative real numbers. Finally, we use the notation \( x(t) = o(t) \) or \( x(t) = O(t) \), where \( x(t) \in \mathbb{R}^n \) and \( t > 0 \), to indicate that \( \frac{\|x(t)\|}{t} \) approaches 0 or a non-zero constant as \( t \to 0 \), respectively.

Next, we review the concepts from variational analysis and generalized differentiation that are commonly used in the following discussion. For further details and references, see \cite{BA00, Mor2018, RWets98}.
Given a set $S\subset X$ and a point $x\in S$, we define the tangent cone of S at $x$
\[
\mathcal{T}_{S}(x):=\{d\in X|~\exists x^k\rightarrow x, t_k\downarrow 0~\mbox{such~~that}~x^k\in S~\mbox{and}~(x^k-x)/t_k\rightarrow d\}.
\]
If the $S\subset X$ is convex, we denote the normal cone by
\[
\mathcal{N}_{S}(x):=\{v\in X|~\langle v,y-x\rangle\leq 0,\forall y\in S\}.
\]
By virtue of \cite[Proposition 3.17]{RWets98}, for a convex set $S$,
\begin{align}\label{2.1}
	y\in \mathcal{N}_{S}(x) \quad\mbox{if~and~only~if}\quad \Pi_{S}(y+x)=x.
\end{align}
Let $\mathcal{F}:X\rightrightarrows W$ be a given multifunction with its  graph given by
$\mbox{gph}\mathcal{F}:=\{(x,w)\in X\times W~|~w\in\mathcal{F}(x)\}$.
Consider an arbitrary $(\bar{x},\bar{w})\in gph \mathcal{F}$ such that $\mathcal{F}$ is locally closed at $(\bar{x},\bar{w})$. The multifunction $\mathcal{F}$ is said to be calm at $\bar{x}$ for $\bar{w}$  if there exists $\kappa>0$ along with $\varepsilon>0$ and $\delta>0$ such that for all $x\in B_{\varepsilon}(\bar{x})$,
\[
\mathcal{F}(x)\cap B_{\delta}(\bar{w})\subset \mathcal{F}(\bar{x})+\kappa\|x-\bar{x}\|\mathbb{B}_{W}.
\]
Equivalently, $\mathcal{F}$ is calm at $\bar{x}$ for $\bar{w}$ if there exists $\kappa>0$, $\varepsilon>0$ and $\delta>0$ such that for all $w\in B_{\delta}(\bar{w})$,
\[
\mbox{dist}(w,\mathcal{F}(\bar{x}))\leq \kappa\mbox{dist}(\bar{x},\mathcal{F}^{-1}(w)\cap B_{\varepsilon}(\bar{x})),	
\]
i.e., the inverse map $\mathcal{F}^{-1}$ is metrically subregular at $\bar{w}$  for $\bar{x}$.

In this paper, we will refer to an important concept from \cite{BA00}.
\begin{definition}\label{d2.2}
	We define  $K$ as being $\mathcal{C}^2-$cone reducible at $y_0\in K$ if there exist a pointed closed convex cone $C\subset H$ in a finite-dimensional space H, a neighborhood $N$ of $y_0$, and a twice continuously differentiable map $\Theta:N\rightarrow H$ such that $\Theta(y_0)=0$, $\Theta'(y_0)$ is surjective, and $K\cap N=\Theta^{-1}(C)\cap N$.  $K$ is said to be  $\mathcal{C}^2-$cone reducible if the above holds for all points $y_0\in K$.
\end{definition}
The $\mathcal{C}^2$-cone reducibility of $K$ is essential for both first and second-order optimality conditions. As mentioned in the introduction, several classes of constraints are known to exhibit
$\mathcal{C}^2$-reducibility.

\subsection{Optimality conditions for conic programmings }

Let $L:X\times H\rightarrow  \mathbb{R}$ denote the Lagrange function associated with problem (\ref{1.1}):
\[
L(x,\lambda):=f(x)+\langle \lambda,g(x)\rangle.
\]The Karush-Kuhn-Tucker (KKT) point $(\bar{x},\bar{\lambda})$ for (\ref{1.1}) satisfies
\begin{align}\label{KKT-system}
	 \nabla _{x}L(\bar{x},\bar{\lambda})=\nabla f(\bar{x})+\nabla g(\bar{x})\bar{\lambda}=0 ,\quad \bar{\lambda}\in \mathcal{N}_{K}(g(\bar{x})).
\end{align}
We refer to $\bar{x}$ as a stationary point if there exists a multiplier $\bar{\lambda}$ such that $(\bar{x},\bar{\lambda})$ satisfies the KKT conditions. The set of such multipliers is denoted by
 $\Lambda(\bar{x})$, defined as:
\begin{align}\label{multiplier set}
\Lambda(\bar{x}):=\{\lambda\in H~|~	\nabla _{x}L(\bar{x},\lambda)=0,\quad \lambda\in \mathcal{N}_{K}(g(\bar{x}))\}.
\end{align}
Since
\begin{align}\label{2.3}
	\nabla_{(x,\lambda)} \mathcal{L}(x, \lambda, c) = \begin{pmatrix}
		\nabla_{x} \mathcal{L}_{c}(x, \lambda) \\
		\nabla_{\lambda} \mathcal{L}_{c}(x, \lambda)
	\end{pmatrix}
	= \begin{pmatrix}
		\nabla f(x) + c \nabla g(x)\left[ g(x) + c^{-1} \lambda - \Pi_{K} \left( g(x) + c^{-1} \lambda \right) \right] \\
		g(x) - \Pi_{K} \left( g(x) + c^{-1} \lambda \right)
	\end{pmatrix},
\end{align}
we have from (\ref{2.1}) that \( (\bar{x}, \bar{\lambda}) \) is a KKT point of (\ref{1.1}) if and only if, for any \( c > 0 \), \( \nabla_{(x, \lambda)} \mathcal{L} (\bar{x}, \bar{\lambda}, c) = 0 \), where \( \mathcal{L} (x, \lambda, c) \) is the augmented Lagrangian function defined in (\ref{augmented function}).

The connection between the optimization problem (\ref{1.1}) and its KKT conditions is well-established \cite{BA00}: if \( \bar{x} \) is a local solution and a suitable constraint qualification holds, then \( (\bar{x}, \bar{\lambda}) \) is a KKT point. Specifically, in convex problems, the KKT conditions are both necessary and sufficient for local and global optimality. In the absence of convexity, local optimality can still be ensured if a second-order sufficient condition is met. For the theory in this paper, we require a second-order condition that depends on the \( \mathcal{C}^2 \)-cone reducibility of the set \( K \). The support function of a closed convex set \( C \) is defined as:
\[
\sigma(y,C) := \sup \{ \langle y, z \rangle \mid z \in C \}.
\]
For a point \( y \in K \) and a direction \( h \in H \), the second-order tangent set to \( K \) in the direction \( h \) is given by:
\[
\mathcal{T}_K^2(y,h) := \{ w \in H \mid \text{dist}(y + t h + \frac{1}{2} t^2 w, K) = o(t^2), \ t \geq 0 \}.
\]
Different definitions of second-order tangent sets exist, such as inner and outer tangent sets. However, for \( \mathcal{C}^2 \)-cone reducible sets \( K \), these definitions coincide \cite[Proposition 3.136]{BA00}, so the distinction is not essential for our purposes.

Let $\bar{x}$ be a stationary point of (\ref{1.1}). Assume $\Lambda(\bar{x})$  is nonempty with $\bar{\lambda}\in \Lambda(\bar{x})$, the critical cone of (\ref{1.1}) is given by
\[
C(\bar{x}):=\{d\in X~|~\mathcal{J}g(\bar{x})d\in\mathcal{T}_{K}(g(\bar{x})),\quad \langle \mathcal{J}g(\bar{x})d,\bar{\lambda} \rangle=0\}.
\]
We now present the second-order sufficient condition for the problem in (\ref{1.1}), which can be derived from \cite[Proposition 7.3]{Mohammadi2021} and \cite[(3.16)]{MS2019}.

\begin{definition}\label{d2.5}
	Let $\bar{x}$ be a stationary point of (\ref{1.1}) and $\bar{\lambda}\in \Lambda(\bar{x})$.  We say the second-order sufficient condition (SOSC) holds at $(\bar{x},\bar{\lambda})$ if
	\begin{align}\label{SOSC}
		\langle\nabla_{xx}^2L(\bar{x},\bar{\lambda})d,d\rangle-\sigma\left(\bar{\lambda},\mathcal{T}^2_{K}(g(\bar{x}),\mathcal{J}g(\bar{x})d)\right)>0,
	\end{align}
	for any     $d\in C(\bar{x})\backslash\{0\}$.
\end{definition}

Note that SOSC in (\ref{SOSC}) is a stricter version than the conventional second-order sufficient condition for (\ref{1.1}). The latter requires that the supremum of the quadratic term in (\ref{SOSC}) over all Lagrange multipliers in \( \Lambda(\bar{x}) \) be positive. It is well-known that this stronger SOSC ensures the quadratic growth of the objective function, meaning there exists a constant \( \kappa > 0 \) such that
\[
f(x) \geq f(\bar{x}) + \kappa \|x - \bar{x}\|^2
\]
for all feasible points \( x \) arround \( \bar{x} \), as shown in \cite[Theorem 3.86]{BA00}.

Also under this SOSC (\ref{SOSC}), we can obtain the following uniform quadratic growth condition of the augmented Lagrangian function.
\begin{theorem}\label{uniform quadric growth}
	Let \( \bar{x} \) be a stationary point of (\ref{1.1}) and \( \bar{\lambda} \in \Lambda(\bar{x}) \). Suppose that the SOSC in (\ref{SOSC}) holds at \( (\bar{x}, \bar{\lambda}) \). Then, there exist positive constants \( \bar{c} \), \( l \), and \( \varepsilon \) such that the uniform quadratic growth condition
	\begin{align}\label{quadratic growth}
	\mathcal{L}(x,\lambda,c)\geq\mathcal{L}(\bar{x},\lambda,c)+\displaystyle\frac{\kappa}{2}\|x-\bar{x}\|^2
\end{align}	
	is satisfied for all \( x \in B_{\varepsilon}(\bar{x}) \), \( \lambda \in \Lambda(\bar{x}) \cap B_{\varepsilon}(\bar{\lambda}) \), \( c \geq \bar{c} \), and \( \kappa \in [0, l) \).
	
\end{theorem}
{\bf{Proof.}}\  It follows from \cite[Theorem 6.2]{Mohammadi2021} that a
$\mathcal{C}^2$-cone reducible set must satisfy parabolic regularity (as defined in \cite[Definition 3.1]{Mohammadi2021}). Combining this result with \cite[Theorem 3.3(iii), Theorem 7.1]{Mohammadi2021} and Remark 7.2 in \cite{Mohammadi2021}, we conclude that the SOSC in (\ref{SOSC}) is equivalent to the existence of a constant \( l' > 0 \) such that the following inequality holds:
\[
\langle\nabla _{xx}^2L(\bar{x},\bar{\lambda})d,d\rangle-\sigma(\bar{\lambda},\mathcal{T}^2_{K}(g(\bar{x}),\mathcal{J}g(\bar{x})d))\geq l'\|d\|^2
\]
for all $d\in C(\bar{x})$. Finally, by \cite[Theorem 4.7]{Hang2023}, for any $l\in  [0,l')$ the conclusion follows. \hfill\textbf{$\Box$}
\subsection{Calmness and semi-isolated calmness of the KKT solution mapping}
We now introduce the perturbed KKT solution mapping $\mathcal{S}_{KKT}:X\times H\rightarrow X\times H$ for problem (\ref{1.1}):
\begin{align}\label{S_KKT}
	\mathcal{S}_{KKT}(a,b):=\{(x,\lambda)\in X\times H~|~a= \nabla_{x}L(x,\lambda),~b\in -g(x)+\mathcal{N}_{K}^{-1}(\lambda)\}.
\end{align}
The form in (\ref{S_KKT}) for perturbed KKT solution mapping has been widely used in discussing the stability of KKT system for conic programming. For example, see  \cite{DingSunZhang2017,LiuPan2019}. Referring to the definitions of calmness for set-valued mappings in the previous section, the calmness condition for $\mathcal{S}_{KKT}$ at the origin for a KKT point $(\bar{x},\bar{\lambda})$ implies there exist  constants $\varepsilon>0$, $\delta>0$ and $l>0$ such that for any $(x,\lambda)\in B_{\varepsilon}(\bar{x},\bar{\lambda})$ and $(a,b)\in X\times H$ satisfying $\|(a,b)\|\leq \delta$,
\begin{align*}
	\mathcal{S}_{KKT}(a,b)\cap B_{\varepsilon}(\bar{x},\bar{\lambda})\subset \mathcal{S}_{KKT}(0,0)+l(\|a\|+\|b\|)\mathbb{B}_{X\times H}.
\end{align*}

Next, we introduce the definition of semi-isolated calmness for $\mathcal{S}_{KKT}$. The concept of semi-isolated calmness of the KKT solution mapping was first formally introduced in \cite[Theorem 4.1]{MS2018} to establish the characterization of noncritical multipliers for the polyhedral problems.
\begin{definition}
	Let $(\bar{x},\bar{\lambda})$ be a KKT point for problem (\ref{1.1}). The multifunction $\mathcal{S}_{KKT}$ is said to be semi-isolated calm  at the origin for $(\bar{x},\bar{\lambda})$ if there exist  constants $\varepsilon>0$, $\delta>0$ and $l>0$ such that for any  $(a,b)\in X\times H$ satisfying $\|(a,b)\|\leq \delta$, we have
	\begin{align*}
		\mathcal{S}_{KKT}(a,b)\cap B_{\varepsilon}(\bar{x},\bar{\lambda})\subset \{\bar{x}\}\times\Lambda(\bar{x})+l(\|a\|+\|b\|)\mathbb{B}_{X\times H}.
	\end{align*}
\end{definition}
From the definitions, it follows that the semi-isolated calmness of the mapping \( \mathcal{S}_{KKT} \) at the origin for \( (\bar{x}, \bar{\lambda}) \) implies its calmness at the origin for \( (\bar{x}, \bar{\lambda}) \).

The residual function \( r : X \times H \to \mathbb{R}_+ \) associated with the KKT system (\ref{KKT-system}) is defined by
\begin{align}\label{residual}
	r(x,\lambda):=\|\nabla_xL(x,\lambda)\|+ \| g(x) - \Pi_K \left( g(x) + \lambda \right)\| .
\end{align}
From (\ref{2.1}), it is straightforward to verify that \( (\bar{x}, \bar{\lambda}) \) is a KKT point if and only if \( r(\bar{x}, \bar{\lambda}) = 0 \). This, combined with the Lipschitz continuity of \( r \) with respect to \( (x, \lambda )\) around \( (\bar{x}, \bar{\lambda}) \), implies the existence of constants \( \varepsilon > 0 \) and \( \hat{\kappa} > 0 \) such that for all \( (x, \lambda) \in B_{\varepsilon}(\bar{x}, \bar{\lambda}) \),
	\begin{align}\label{r}
		r(x,\lambda)\leq \hat{\kappa} \left( \| x - \bar{x} \| + \mbox{dist}(\lambda,\Lambda(\bar{x})) \right),
	\end{align}
where \( \Lambda(\bar{x}) \) is defined in (\ref{multiplier set}).

Below, we show that the opposite inequality in (\ref{r}), which is  equivalent to the semi-isolated calmness of the KKT solution mapping
$\mathcal{S}_{KKT}$. The result is established in \cite[Theorem 3.1]{LiuPan2019}.
\begin{proposition}\label{semi-isolated}
	Let \( (\bar{x}, \bar{\lambda}) \) be a KKT point of (\ref{1.1}).  $\mathcal{S}_{KKT}$ is semi-isolated calm at the origin for $(\bar{x},\bar{\lambda})$ if and only if there are $\varepsilon>0$ and $\tau>0$ such that the estimate
	\begin{align*}
		\|x-\bar{x}\|+\mbox{dist}\left(\lambda,\Lambda(\bar{x})\right) \leq \tau r(x,\lambda)
	\end{align*}
	is satisfied for all pairs $(x,\lambda)\in B_{\varepsilon}(\bar{x},\bar{\lambda})$.
\end{proposition}
For any fixed $\bar{x}$, the multiplier multifunction $M_{\bar{x}}:X\times H\rightarrow H$ is defined as:
\begin{align}\label{M_KKT}
	M_{\bar{x}}(a,b):=\{\lambda\in H~|~~a= \nabla_{x}L(\bar{x},\lambda),~\lambda\in \mathcal{N}_{K}(g(\bar{x})-b)\}.
\end{align}
For a KKT point \( (\bar{x}, \bar{\lambda}) \), the calmness condition for \( M_{\bar{x}} \) at the origin with respect to \( \bar{\lambda} \) requires the existence of positive constants \( l \), \( \varepsilon \), and \( \delta \) such that
\[
M_{\bar{x}}(a,b)\cap B_{\varepsilon}(\bar{\lambda})\subset \Lambda(\bar{x})+l(\|a\|+\|b\|)\mathbb{B}_{H},
\]
whenever $\|(a,b)\|\leq \delta$. It follows from \cite[Theorem 5.9]{MS2019} that the semi-isolated calmness of the solution map $\mathcal{S}_{KKT}$ at origin for $(\bar{x},\bar{\lambda})$ can be achieved if both  $M_{\bar{x}}(\cdot,\cdot)$  is calm at the origin for $\bar{\lambda}$ and the SOSC in (\ref{SOSC})  holds at $(\bar{x},\bar{\lambda})$.
\begin{proposition}
Assume \( (\bar{x}, \bar{\lambda}) \) to be a KKT point of the problem in (\ref{1.1}).
Suppose that the multiplier mapping $M_{\bar{x}}(\cdot,\cdot)$ defined in (\ref{M_KKT}) is calm at the origin for $\bar{\lambda}$ and the  SOSC in (\ref{SOSC}) holds at $(\bar{x},\bar{\lambda})$. Then the solution map $\mathcal{S}_{KKT}$ is semi-isolated calm at origin for $(\bar{x},\bar{\lambda})$.
\end{proposition}

\subsection{An iterative framework for generalized equation}
We will now provide a brief overview of an abstract iterative framework introduced by Fischer in \cite{Fischer2002}, designed to solve generalized equations with nonisolated solutions. This framework plays a vital role in the convergence analysis presented in this paper, as the convergence results are proven based on this framework. Let \( H \) and \( H' \) be two finite-dimensional Hilbert spaces. Given the single-valued mapping \( \Phi : H \to H' \) and set-valued mapping \( G : H \rightrightarrows H' \),  the generalized equation we consider is formulated as follows:
\begin{align}\label{general equation}
0 \in \Phi(u) + G(u).
\end{align}
Let \( \Gamma : H' \rightrightarrows H \) be the solution map for the canonical perturbation of (\ref{general equation}) as follows:
\begin{align}\label{perturbation}
	\Gamma(s) = \left\{ u \in H \mid s \in \Phi(u) + G(u) \right\}, \quad s \in H'.
\end{align}
In iterative algorithms designed to solve (\ref{general equation}), a sequence \( \{u^k\} \) is often generated by solving subproblems:
\begin{align}\label{approximal problem}
	0 \in \mathcal{A}(u, u^k, p_k) + G(u),
\end{align}
where \(\mathcal{A} : H \times H \times P \rightrightarrows H' \) is a set-valued mapping.
 In this approach, the single-valued component of (\ref{general equation}) is approximated by the set-valued mapping \( \mathcal{A}\), and approximates \( \Phi \) near the current iterate \( u^k \), while the set-valued component remains the same.  To guarantee convergence and determine the rate of convergence for the sequence \( \{u^k\} \), it is essential to select a solution to (\ref{approximal problem}) that possesses specific properties. In particular, for a given iterate \( u^k \) and parameter \( p_k \), we choose the next iterate \( u^{k+1} \) to be sufficiently near \( u^k \) so that:
\begin{align}\label{iterative system}
 0 \in \mathcal{A}(u^{k+1}, u^k, p_k) + G(u^{k+1})  \quad \text{and} \quad \| u^{k+1} - u^k \| \leq \alpha \, \text{dist}(u^k, \Gamma(0)),
\end{align}
where \( \alpha > 0 \) is a constant, and \( \Gamma(0) \) is the solution set to (\ref{general equation}).

The following result, established in \cite[Theorem 2.1]{Hang2023}, provides the analysis of the local convergence for the sequence\( \{u^k\} \) under relatively mild assumptions. This theorem was originally presented in \cite[Theorem 4.1]{IK2013}. In that work, assumption (a) was stated in a slightly stronger form. Specifically, it was assumed that the solution mapping $\Gamma$ is upper Lipschitzian, rather than calm. Later, the upper Lipschitzian condition was substituted with calmness in \cite[Theorem 7.13]{IS2014},
\begin{theorem}
Suppose that \( \Gamma(0) \) is locally closed around \( \bar{u} \), where \( \Gamma \) is defined by (\ref{general equation}). Let \( \bar{u} \in \Gamma(0) \). Then, there exist a constant \( \alpha > 0 \) such that  the following statements hold:\\
	(a) The solution mapping \( \Gamma \) is calm at \( 0 \in H' \)  for  \( \bar{u} \) with a constant \( \ell_1 > 0 \);\\
	(b) There exists a constant \( \varepsilon_1 > 0 \) such that for every \( \tilde{u} \in B_{\varepsilon_1}(\bar{u}) \) and for any \( p \in P \), the set
	\[
	\{ u \in H \mid 0 \in \mathcal{A}(u, \tilde{u}, p) + G(u)\}\cap\{u\in H \mid \| u - \tilde{u} \| \leq \alpha \, \text{dist}(\tilde{u}, \Gamma(0)) \}
	\]
	is non-empty;\\
	(c) There exist constants \( \varepsilon_2 > 0 \) and a function \( \omega : H \times H \times P \to \mathbb{R}_+ \) such that
	\[
	\sup\{ \omega(u, \tilde{u}, p) \mid \tilde{u} \in B_{\varepsilon_2}(\bar{u}), \, \| u - \tilde{u} \| \leq \alpha \, \text{dist}(\tilde{u}, \Gamma(0)), p \in P \} < \frac{1}{\ell_1},
	\]
	where \( \ell_1 > 0 \) is the constant from condition (a), and the following estimate holds:
	\[
	\sup\{\|w\| \mid w \in \Phi(u) - \mathcal{A}(u, \tilde{u}, p) \} \leq \omega(u, \tilde{u}, p) \, \text{dist}(\tilde{u}, \Gamma(0))
	\]
	for all \( \tilde{u} \in B_{\varepsilon_2}(\bar{u}) \), all \( u \in H \) with \( \| u - \tilde{u} \| \leq \alpha \, \text{dist}(\tilde{u}, \Gamma(0)) \), and all \( p \in P \).
	
	Then, there exists a constant \( \varepsilon_0 > 0 \) such that for any initial point \( u^0 \in B_{\varepsilon_0}(\bar{u}) \) and any sequence \( \{ p_k \} \subset P \), the sequence \( \{ u_k \} \) generated by the iterative scheme defined by equation (\ref{iterative system}) converges to some \( \hat{u} \in \Gamma(0) \). The convergence of the sequence \( \{ u_k \} \) to \( \hat{u} \) and the process in which the distance \( \{ \text{dist}( u_k, \Gamma(0) ) \} \) approaches zero both exhibit Q-linear convergence. If, as \( k \to \infty \), it holds that \( \omega(u_{k+1}, u_k, p_k) \to 0 \), then Q-superlinear convergence is further achieved.
\end{theorem}

The main convergence result of this paper is derived from the above theorem. Therefore, the primary task in the subsequent sections is to verify that, under suitable conditions, all three assumptions in Theorem 2.2 can be satisfied.
The final goal of this section is to establish that the iterate sequence \(\{(x^k, \lambda^k)\}\), generated by the inexact PALM algorithm, fits the iterative pattern discussed earlier. To achieve this, we define two mappings: \(\Phi : X \times H  \to X \times H \) and \(G : X \times H  \rightrightarrows X \times H \), as follows:
\begin{align}\label{GE2}
	\Phi(x,\lambda)=\begin{pmatrix}
		\nabla_x L(x, \lambda)\\-g(x)
	\end{pmatrix}\quad \mbox{and}\quad G(x,\lambda)=\begin{pmatrix}
		0\\\mathcal{N}_{K}^{-1}(\lambda)
	\end{pmatrix}.
\end{align}
These mappings clearly align with the generalized equation described in (\ref{general equation}). We can verify that \((\bar{x}, \bar{\lambda})\) is a solution to the generalized equation
\begin{align}\label{GE3}
	0\in \Phi(x,\lambda)+G(x,\lambda)
\end{align} if and  only if \((\bar{x}, \bar{\lambda})\) is a KKT point of (\ref{1.1}), and $\mathcal{S}_{KKT}$ is the solution mapping for the canonical perturbation of (\ref{GE3}).
Moreover, the inexact primal update \(x^{k+1}\), as defined by (\ref{x-pupdate}), can be equivalently expressed through the following equation:
\[
0 \in \nabla_x \mathcal{L}(x^{k+1}, \lambda^k, c_k) + c_k^{-1}(x^{k+1} - x^k) + \epsilon_k \mathbb{B}_X,
\]
By (\ref{2.3}), this is equivlent to
\[
0 \in \nabla f(x^{k+1}) + \nabla g(x^{k+1})\left[ c_k g(x^{k+1}) + \lambda^k - \Pi_K\left( c_k g(x^{k+1}) + \lambda^k \right) \right] + c_k^{-1}(x^{k+1} - x^k) + \epsilon_k \mathbb{B}_X,
\]
which simplifies to
\[
0 \in  \nabla_x L(x^{k+1}, \lambda^{k+1}) + c_k^{-1}(x^{k+1} - x^k) + \epsilon_k \mathbb{B}_X.
\]

Next, to show that the iterate \((x^{k+1}, \lambda^{k+1})\) satisfies  (\ref{approximal problem}), we define a mapping \(\mathcal{A}: X \times H\times X \times H \times \mathbb{R}_+ \times (0, +\infty) \rightrightarrows X \times H \times X\) as
\begin{align}\label{A1}
	\mathcal{A} (x, \lambda, \tilde{x}, \tilde{\lambda}, \epsilon, c) :=\begin{pmatrix}
		\nabla_x L(x, \lambda) + c^{-1}(x - \tilde{x}) + \epsilon \mathbb{B}_X\\
		- g(x) + c^{-1}(\lambda - \tilde{\lambda})
	\end{pmatrix},
\end{align}
where \((\tilde{x}, \tilde{\lambda}, \epsilon, c) \in X \times H \times X \times \mathbb{R}_+ \times (0, +\infty)\) is a given point.  Then the iterate \((x^{k+1}, \lambda^{k+1})\) satisfies the subproblem
\[
0 \in \mathcal{A}(x, \lambda, x^k, \lambda^k, \epsilon_k, c_k) + G(x, \lambda).
\]

\section{Local convergence analysis for PALM}
 \setcounter{equation}{0}
This section primarily aims to analyze the local convergence of the inexact PALM algorithm outlined below.

\textbf{Algorithm 4.1 (PALM)}.
Given an starting point \((x^0, \lambda^0) \in X \times H\), a constant \(\bar{c} > 0\), and a constant \(\alpha > 0\), choose a sequence \(\{c_k\}\) with \(c_k \geq \bar{c}\) for all \(k\),  along with  a function \(\epsilon : \mathbb{R}_+ \to \mathbb{R}_+\), which satisfied \(\epsilon(t) = o(t)\). Set \(k := 0\). Then, perform the following steps:
\begin{enumerate}
	\item If \((x^k, \lambda^k)\) satisfies a appropriate termination criterion, stop.
	\item Otherwise, set \(\epsilon_k := \epsilon(r(x^k, \lambda^k))\), where \(r\) is defined in  (\ref{residual}), and update  \((x^{k+1}, \lambda^{k+1})\) according to (\ref{x-pupdate}) and (\ref{lambda-update}), such that:
\begin{align}\label{esta}
	 \|x^{k+1} - x^k\|+ \|\lambda^{k+1} - \lambda^k\|\leq \alpha r(x^k, \lambda^k).
\end{align}
	\item Set \(k \leftarrow k + 1\) and return to Step 1.
\end{enumerate}

The approach to guaranteeing the convergence of the sequence \(\{(x^k, \lambda^k)\}\), generated by the PALM algorithm in Algorithm 4.1, is outlined in Theorem 2.2. This section is dedicated to verifying assumption (b) of Theorem 2.2, which involves two key tasks: first, we verify the solvability of the subproblem in (\ref{x-update}), and second, we obtain an error bound estimate for the iterates of the PALM algorithm. We begin by addressing the solvability of the subproblems. While the solvability of the ALM subproblems for the piecewise linear-quadratic composite optimization problem and NLSDP was established in \cite[Proposition 5.2]{Hang2021} and \cite[Proposition 4]{Wang2024}, respectively, the solvability of PALM subproblems has not yet been studied. In this work, we provide a detailed proof of the solvability of the subproblems in the context of the PALM method.

%
\begin{proposition}\label{S-isolated-calm}
Let \( \bar{x} \) be a stationary point of (\ref{1.1}) and \( \bar{\lambda} \in \Lambda(\bar{x}) \). Assume that the SOSC in (\ref{SOSC}) holds at \( (\bar{x}, \bar{\lambda}) \). Let \( \bar{c} \) and \( \varepsilon \) be positive constants from Theorem \ref{uniform quadric growth}. Then, there exist positive constants \( \hat{l} \) and \( \hat{\varepsilon} \in (0, \varepsilon) \) such that for any \( c \geq \bar{c} \), the local optimal solution mapping \( \mathcal{S}_{c}: H \times X \rightrightarrows X \) defined by
\begin{align}\label{subproblem}
\mathcal{S}_{c}(\lambda, v) := \arg \min_{ x \in B_{\hat{\varepsilon}}(\bar{x})} \left\{ \mathcal{L}(x, \lambda, c) + \frac{1}{2c} \|x - v\|^2 \right\}, \quad \lambda \in H, \quad v \in X,
\end{align}
enjoys the uniform isolated calmness property
\begin{align}\label{isolated calm for S}
	\mathcal{S}_{c}(\lambda,v) \subset \{ \bar{x} \} +\hat{l} \left(\|\lambda - \bar{\lambda}\|+\|v -\bar{x}\|\right)\mathbb{B}_{X}
\end{align}
and satisfies the inclusion
\[
\emptyset \neq \mathcal{S}_{c}(\lambda, v) \subset \text{int} B_{\hat{\varepsilon}}(\bar{x}),
\]
for all \( \lambda \in B_{\hat{\varepsilon}/(3\hat{l})}(\bar{\lambda}) \) and \( v \in B_{\hat{\varepsilon}/(3\hat{l})}(\bar{x}) \).
\end{proposition}
\textbf{Proof.}\ Let $\varepsilon >0$ from Theorem \ref{uniform quadric growth}, since \(g(x)\) is twice differentiable at \(\bar{x}\), there are constants $\hat{\varepsilon} \in (0, \varepsilon)$ and \(\tau_g> 0\) such that
\begin{align}\label{g-lip}
\|g(x) - g(\bar{x})\| \leq \tau_g \| x - \bar{x} \| \quad \text{for all} \quad x \in B_{\hat{\varepsilon}}(\bar{x}).
\end{align}
It follows from Theorem 2.1 that the  SOSC in (\ref{SOSC}) at  $(\bar{x},\bar{\lambda})$ implies the the quadratic growth condition (\ref{quadratic growth}), which means  \(\mathcal{S}_{c}(\bar{\lambda},\bar{x}) = \{ \bar{x} \}\) for all \(c \geq \bar{c}\). Therefore, by applying \cite[Theorem 4.16]{Rudin1976}, it can be deduced from the continuity of \( \mathcal{L} \) and the compactness of \( B_{\hat{\varepsilon}}(\bar{x}) \) that \( \mathcal{S}_{c}(\lambda, v) \) is non-empty.

 We know from \cite{RWets98} that the augmented Lagrangian $\mathcal{L}(x, \lambda, c)$ is concave with respect to $\lambda$.
  Fix any $u \in \mathcal{S}_{c}(\lambda,v)$, combined with (\ref{2.3}), there exist constants $\bar{c}>0$ and $l>0$ from Theorem \ref{uniform quadric growth} such that for all $c \geq \bar{c}$, $\lambda \in H$, $v\in X$ and $\kappa\in [0,l)$,
\begin{align*}
\mathcal{L}(u, \lambda, c) +\displaystyle\frac{1}{2c}\|u-v\|^2
&\geq \mathcal{L}(u, \bar{\lambda}, c) - \langle\nabla_{\lambda} \mathcal{L}(u, \lambda, c), \bar{\lambda} - \lambda \rangle+\displaystyle\frac{1}{2c}\|u-v\|^2\\
&= \mathcal{L}(u, \bar{\lambda}, c)- \langle g(u) - \Pi_K \left( g(u) + c^{-1} \lambda \right),\bar{\lambda} - \lambda\rangle+\displaystyle\frac{1}{2c}\|u-v\|^2\\
&\geq f(\bar{x}) + \displaystyle\frac{\kappa}{2}\| u - \bar{x} \|^2 - \langle g(u) - \Pi_K \left( g(u) + c^{-1} \lambda \right),\bar{\lambda} - \lambda\rangle+\displaystyle\frac{1}{2c}\|u-v\|^2.
\end{align*}
On the other hand, by the optimality of $u$, we have:
\begin{align*}
\mathcal{L}(u, \lambda, c)+\displaystyle\frac{1}{2c}\|u-v\|^2 &\leq \mathcal{L}(\bar{x}, \lambda, c)+\displaystyle\frac{1}{2c}\|\bar{x}-v\|^2\\
& = f(\bar{x}) + \frac{c}{2} \, \text{dist}^2 \left( g(\bar{x}) + c^{-1}\lambda, K \right) - \frac{1}{2c}  \|\lambda \|^2 +\displaystyle\frac{1}{2c}\|\bar{x}-v\|^2\\
&\leq f(\bar{x})+\displaystyle\frac{1}{2c}\|\bar{x}-v\|^2.
\end{align*}
Combining the above two inequalities, we get:
\[
\displaystyle\frac{\kappa}{2}\| u - \bar{x} \|^2 - \langle g(u) - \Pi_K \left( g(u) + c^{-1} \lambda \right),\bar{\lambda} - \lambda\rangle+\displaystyle\frac{1}{2c}\|u-v\|^2\leq \displaystyle\frac{1}{2c}\|\bar{x}-v\|^2,
\]
hence
\begin{align}\label{relation}
\displaystyle\frac{\kappa}{2}\| u - \bar{x} \|^2 +\displaystyle\frac{1}{2c}\|u-\bar{x}\|^2\leq \langle g(u) - \Pi_K \left( g(u) + c^{-1} \lambda \right),\bar{\lambda} - \lambda\rangle+\displaystyle\frac{1}{c} \langle u-\bar{x},v-\bar{x}\rangle.
\end{align}
It follows from (\ref{2.1}) that \( \bar{\lambda} \in \mathcal{N}_K (g(\bar{x})) \) implies
$g(\bar{x}) = \Pi_K \left( g(\bar{x}) + c^{-1} \bar{\lambda} \right)$.
Then
\begin{align*}
\|g(u) - \Pi_K \left( g(u) + c^{-1} \lambda \right) \|&= \|g(u) - g(\bar{x}) + \Pi_K \left( g(\bar{x}) + c^{-1} \bar{\lambda } \right) - \Pi_K \left( g(u) + c^{-1} \lambda  \right)\|\\
&\leq 2 \| g(u) - g(\bar{x}) \| + c^{-1} \|\lambda-\bar{\lambda}  \|\\
&\leq 2 \tau_g \| u - \bar{x} \| + c^{-1} \| \lambda-\bar{\lambda} \|,
\end{align*}
where the last inequality comes from (\ref{g-lip}). Using this and (\ref{relation}) tells us that
\begin{align*}
	\displaystyle\frac{\kappa}{2}\| u - \bar{x} \|^2 +\displaystyle\frac{1}{2c}\|u-\bar{x}\|^2\leq &\| g(u) - \Pi_K \left( g(u) + c^{-1} \lambda \right)\|\|\lambda-\bar{\lambda} \|+\displaystyle\frac{1}{c} \| u-\bar{x}\|\|v-\bar{x}\|\\
	\leq& 2 \tau_g \| u - \bar{x} \|\| \lambda-\bar{\lambda}  \| + \displaystyle\frac{1}{c}  \| \lambda-\bar{\lambda}  \|^2+\displaystyle\frac{1}{c} \| u-\bar{x}\|\|v-\bar{x}\|,
\end{align*}
which can be written in the equivalent form as
\[(c\kappa+1)\| u - \bar{x} \|^2-4c\tau_g\| u - \bar{x} \|\| \lambda-\bar{\lambda} \|-2 \| u-\bar{x}\|\|v-\bar{x}\|\leq 2\| \lambda-\bar{\lambda}  \|^2.
\]
By completing the square in the equation above, we arrive at the following expression:
\begin{align}\label{sqrt form}
c\kappa\left[\| u - \bar{x} \|-\displaystyle\frac{2\tau_g}{\kappa}\| \lambda-\bar{\lambda} \|\right]^2+\left[\| u - \bar{x} \|-\| v-\bar{x} \|\right]^2\leq \left(2+\displaystyle\frac{4c\tau_g^2}{\kappa}\right)\| \lambda-\bar{\lambda} \|^2+\|v-\bar{x}\|^2.
\end{align}
It is well known that for any constant $a$ and $b$, the following triangle inequality holds:
\[\displaystyle\frac{a+b}{\sqrt{2}}\leq\sqrt{a^2+b^2},\]
and for any positive constant $a$ and $b$,  the inequality
\[\sqrt{a^2+b^2}\leq a+b\]
also holds. Using these triangle inequalities, we can derive the following relationships from the left-hand side of inequality  (\ref{sqrt form}):
\begin{align}\label{left}
	\displaystyle\frac{\sqrt{c\kappa}\|u-\bar{x}\|-\displaystyle\frac{2\tau_g\sqrt{c\kappa}}{\kappa}\| \lambda-\bar{\lambda} \|+\| u-\bar{x}\|-\|v-\bar{x}\|}{\sqrt{2}}\leq\sqrt{c\kappa\left[\| u - \bar{x} \|-\displaystyle\frac{2\tau_g}{\kappa}\| \lambda-\bar{\lambda} \|\right]^2+\left[\| u - \bar{x} \|-\| v-\bar{x} \|\right]^2}
\end{align}
Similarly, from the right-hand side of inequality (\ref{sqrt form}), we obtain the following:
\begin{align}\label{right}
	\sqrt{\left(2+\displaystyle\frac{4c\tau_g^2}{\kappa}\right)\|\lambda-\bar{\lambda} \|^2+\|v-\bar{x}\|^2}\leq \sqrt{2+\displaystyle\frac{4c\tau_g^2}{\kappa}}\|\lambda-\bar{\lambda}\|+\|v-\bar{x}\|.
\end{align}
Thus,  inequalities (\ref{sqrt form}), (\ref{left}) and (\ref{right}) in turn give us the estimate
\[
\| u - \bar{x} \|\leq \displaystyle\frac{\sqrt{4+\displaystyle\frac{8c\tau_g^2}{\kappa}}+\displaystyle\frac{2\tau_g\sqrt{c}}{\sqrt{\kappa}}}{\sqrt{c\kappa}+1}\| \lambda-\bar{\lambda}  \|+\displaystyle\frac{\sqrt{2}+1}{\sqrt{c\kappa}+1}\|v-\bar{x}\|.
\]
Let $\hat{l}:=\max\left\{\displaystyle\frac{\sqrt{4+\displaystyle\frac{8c\tau_g^2}{\kappa}}+\displaystyle\frac{2\tau_g\sqrt{c}}{\sqrt{\kappa}}}{\sqrt{c\kappa}+1},\displaystyle\frac{\sqrt{2}+1}{\sqrt{c\kappa}+1}\right\}$,  for any $\lambda\in B_{\hat{\varepsilon}/(3\hat{l})}(\bar{\lambda})$, $v\in B_{\hat{\varepsilon}/(3\hat{l})}(\bar{x})$,
 we have
\[
\| u - \bar{x} \|\leq \hat{l}\| \lambda-\bar{\lambda}  \|+\hat{l}\|v-\bar{x}\|\leq 2\hat{l}\cdot \displaystyle\frac{\hat{\varepsilon}}{3\hat{l}}<\hat{\varepsilon},
\]
which completes the proof. \hfill\textbf{$\Box$}

Next, we focus on the second condition specified in assumption (b) of Theorem 2.2, which involves deriving an error bound estimate for the successive terms of our inexact PALM. The following results draw an analogy to \cite[Propostion 5.4]{Hang2023}, which primarily examines the ALM applied to composite optimization problems.
\begin{proposition}
	Let $\bar{x}$
be a stationary point of (\ref{1.1}), and let  $\bar{\lambda}\in \Lambda(\bar{x})$. Assume that the SOSC in (\ref{SOSC}) holds at
$(\bar{x},\bar{\lambda})$ and $M_{\bar{x}}$ is calm at the origin for $\bar{\lambda}$. Let \( \bar{c} \) be the constant from Theorem \ref{uniform quadric growth}. Then, there exist positive constants \( \varepsilon \) and \( \alpha \) such that for any \( c \geq \bar{c} \), for any  \( (x, \lambda) \in B_{\varepsilon }(\bar{x}, \bar{\lambda}) \) satisfying \( r(x, \lambda) > 0 \), and for any optimal solution \( u\in \mathcal{S}_{c}(\lambda,x) \) of the subproblem (\ref{subproblem}),	the following error bound estimate holds:
\begin{align}\label{error bound}
	\|u - x\| + \|\lambda_u - \lambda\| \leq \alpha r(x, \lambda),
\end{align}
where $\lambda_u:=c g(u)+\lambda - \Pi_K \left( c g(u) +  \lambda \right)$ and the residual function $r(x,\lambda)$ is given by (\ref{residual}).
\end{proposition}
\textbf{Proof.}
Given that the SOSC in (\ref{SOSC}) holds, we can apply Theorem \ref{uniform quadric growth} to determine positive constants \( \bar{c} \), \( \varepsilon_1 \), \( l \), and \( \kappa \in (0, l) \) such that the uniform quadratic growth condition in (\ref{quadratic growth}) holds for all \( x \in B_{\varepsilon_1}(\bar{x}) \), all \( \lambda \in \Lambda(\bar{x}) \cap B_{\varepsilon_1}(\bar{\lambda}) \), and for all \( c \geq \bar{c} \). Based on Proposition \ref{S-isolated-calm}, there exist positive constants \( \hat{l} \) and \( \hat{\varepsilon} \in (0, \varepsilon_1) \) such that \( \mathcal{S}_c \) satisfies the uniform isolated calmness property in (\ref{isolated calm for S}). Furthermore, we have \( \mathcal{S}_c(\lambda, x) \subset \text{int} B_{\hat{\varepsilon}}(\bar{x}) \) for all \( \lambda \in B_{\hat{\varepsilon} / (3\hat{l})} (\bar{\lambda}) \), \( x \in B_{\hat{\varepsilon} / (3\hat{l})} (\bar{x}) \), and for all \( c \geq \bar{c} \). As a result, for every \( \tilde{\lambda} \in B_{\hat{\varepsilon}/(3\hat{l})}(\bar{\lambda}) \), \( \tilde{x} \in B_{\hat{\varepsilon}/(3\hat{l})}(\bar{x}) \), and for all \( c \geq \bar{c} \), any optimal solution \( u \in \mathcal{S}_{c}(\tilde{\lambda}, \tilde{x}) \) must satisfy the first-order optimality condition:
\begin{align}\label{optimial_subproblem}
\nabla_x \mathcal{L}(u, \tilde{\lambda}, c) + c^{-1}(u - \tilde{x }) = 0.
\end{align}

Suppose, for the sake of contradiction, that the inequality in (\ref{error bound}) does not hold. This assumption implies the existence of a sequence \(\{(x^k, \lambda^k, c_k)\} \) in \( X \times H \times [\bar{c}, +\infty)\), where \((x^k, \lambda^k) \to (\bar{x}, \bar{\lambda})\), and an optimal solution \(u^k\in \mathcal{S}_{c_k}(\lambda^k,x^k)\), such that
\begin{align*}
\|u^k - x^k\| + \|c_{k} g(u^k) - \Pi_K \left(c_{k} g(u^k) +  \lambda^k \right)\|  > k r_k,
\end{align*}
where \(r_k := r(x^k, \lambda^k) > 0\) and $r_k\rightarrow 0$ as $k\rightarrow +\infty$.  Denoting by
\[
\beta_k:=\|u^k - x^k\| +\|c_{k} g(u^k) - \Pi_K \left(c_{k} g(u^k) +  \lambda^k \right)\|,
\]
we get
\[\displaystyle\frac{r_k}{\beta_k}\rightarrow 0\quad \text{as} \quad k \to \infty.\]
It follows the definition of the residual function $r$ in (\ref{residual}) that
\begin{align*}
\displaystyle\frac{\nabla_x L(x^k, \lambda^k)}{\beta_k} \rightarrow 0\quad \text{and} \quad  \displaystyle\frac{g(x^k)-\Pi_K \left( g(x^k) + \lambda^k \right)}{\beta_k} \rightarrow 0 \quad \text{as} \quad k \to \infty.
\end{align*}
By passing to a subsequence if needed, we can identify a pair \((\xi, \eta) \in X \times H\) such that
\begin{align}\label{3}
\frac{u^k - x^k}{\beta_k} \to \xi~\text{and} ~\frac{c_{k} g(u^k) - \Pi_K \left(c_{k} g(u^k) +  \lambda^k \right)}{\beta_k} \to \eta ~\text{with} ~(\xi,\eta ) \neq 0.
\end{align}
On the other hand, it follows from Proposition 2.2 that the calmness of $M_{\bar{x}}$ at the origin for $\bar{\lambda}$ and the SOSC in (\ref{SOSC}) ensure the semi-isolated calmness of the solution map $\mathcal{S}_{KKT}$ at the origin for $(\bar{x},\bar{\lambda})$. By Proposition \ref{semi-isolated}, there exist a positive constant $\kappa'$ such that for any $(x^k,\lambda^k)$ sufficiently close to $(\bar{x},\bar{\lambda})$,
\begin{align}\label{4}
\|x^k-\bar{x}\|+\mbox{dist}(\lambda^k, \Lambda(\bar{x}))\leq \kappa'r(x^k,\lambda^k).
\end{align}
Since the set of Lagrange multipliers \(\Lambda(\bar{x})\) is convex and closed, \(\Pi_{\Lambda(\bar{x})}(\lambda^k)\) exists and is a singleton. Let \(\tilde{\lambda}^k := \Pi_{\Lambda(\bar{x})}(\lambda^k)\). From (\ref{4}), we can deduce that \( x^k - \bar{x} = O(r_k) \) and \( \lambda^k - \tilde{\lambda}^k = O(r_k) \) for sufficiently large \(k\). Therefore, we obtain
\begin{align}\label{5}
	\frac{x^k - \bar{x}}{\beta_k} \to 0 \quad \text{and} \quad \frac{\lambda^k - \tilde{\lambda}^k}{\beta_k} \to 0 \quad \text{as} \quad k \to \infty,
\end{align}
which, combined with the fact that \( \lambda^k \to \bar{\lambda} \), implies that \( \tilde{\lambda}^k \to \bar{\lambda} \). Hence, for all sufficiently large \(k\), we have \( \tilde{\lambda}^k \in \Lambda(\bar{x}) \cap B_{\hat{\varepsilon}}(\bar{\lambda}) \).
This, along with the fact that \(u^k\in B_{\hat{\varepsilon}}(\bar{x})\), \(c_k \geq \bar{c}\), (\ref{quadratic growth}), and the concavity of
 \(\lambda \to \mathcal{L}(u^k, \lambda, c_k)\) with respect to $\lambda$, leads to
\begin{align*}
\|u^k - \bar{x}\|^2 &\leq \displaystyle\frac{2}{\kappa} \left( \mathcal{L}(u^k, \tilde{\lambda}^k, c_k) - \mathcal{L}(\bar{x}, \tilde{\lambda}^k, c_k) \right) \\&\leq \displaystyle\frac{2}{\kappa}\left( \mathcal{L}(u^k, \lambda^k, c_k) + \langle\nabla_\lambda \mathcal{L}(u^k, \lambda^k, c_k), \tilde{\lambda}^k - \lambda^k \rangle- \mathcal{L}(\bar{x}, \tilde{\lambda}^k, c_k) \right)\\
&=\displaystyle\frac{2}{\kappa}\left( \mathcal{L}(u^k, \lambda^k, c_k) + \langle g(u^k) - \Pi_K \left( g(u^k) + c_{k}^{-1} \lambda^k \right), \tilde{\lambda}^k - \lambda^k \rangle- \mathcal{L}(\bar{x}, \tilde{\lambda}^k, c_k) \right).
\end{align*}
 Given that
 \(u^k \in \mathcal{S}_{c_k}(\lambda^k,x^k)\), we can conclude that
\begin{align*}
\mathcal{L}(u^k, \lambda^k, c_k)+\displaystyle\frac{1}{2c_k}\|u^k-x ^k\|^2 &\leq \mathcal{L}(\bar{x}, \lambda^k, c_k) +\displaystyle\frac{1}{2c_k}\|\bar{x}-x ^k\|^2 \\
&\leq f(\bar{x}) +\displaystyle\frac{1}{2c_k}\|\bar{x}-x ^k\|^2\\
&= \mathcal{L}(\bar{x}, \tilde{\lambda}^k, c_k)+\displaystyle\frac{1}{2c_k}\|\bar{x}-x ^k\|^2,
\end{align*}
where the equality comes from the definition of augmented Lagrangian $\mathcal{L}$ and \(\tilde{\lambda}^k \in \Lambda(\bar{x})\). By combining the last two inequalities, it holds that
\begin{align*}
\|u^k - \bar{x}\|^2 &\leq \displaystyle\frac{2}{\kappa}\left( \left\langle g(u^k) - \Pi_K \left( g(u^k) + c_{k}^{-1} \lambda^k \right), \tilde{\lambda}^k - \lambda^k \right\rangle + \mathcal{L}(u^k, \lambda^k, c_k) - \mathcal{L}(\bar{x}, \tilde{\lambda}^k, c_k) \right)\\
&\leq \displaystyle\frac{2}{\kappa}\left( \left\langle g(u^k) - \Pi_K \left( g(u^k) + c_{k}^{-1} \lambda^k \right), \tilde{\lambda}^k - \lambda^k \right\rangle +\displaystyle\frac{1}{2c_k}\|\bar{x}-x ^k\|^2-\displaystyle\frac{1}{2c_k}\|u^k-x ^k\|^2\right)\\
&=\displaystyle\frac{2}{\kappa}\left( \left\langle g(u^k) - \Pi_K \left( g(u^k) + c_{k}^{-1} \lambda^k \right), \tilde{\lambda}^k - \lambda^k \right\rangle -\displaystyle\frac{1}{2c_k}\|u^k-\bar{x}\|^2-\displaystyle\frac{1}{c_k}\left\langle u^k-\bar{x},\bar{x}-x ^k\right\rangle\right)\\
&\leq\displaystyle\frac{2}{\kappa}\| g(u^k) - \Pi_K \left( g(u^k) + c_{k}^{-1} \lambda^k \right)\|\|\tilde{\lambda}^k - \lambda^k \| -\displaystyle\frac{1}{\kappa c_k}\|u^k-\bar{x}\|^2+\displaystyle\frac{2}{\kappa c_k}\|u^k-\bar{x}\|\|\bar{x}-x ^k\|.
\end{align*}
Then we have
\[
(\kappa c_k+1)\|u^k - \bar{x}\|^2\leq 2c_k\| g(u^k) - \Pi_K \left( g(u^k) + c_{k}^{-1} \lambda^k \right)\|\|\tilde{\lambda}^k - \lambda^k\|+2\|u^k-\bar{x}\|\|\bar{x}-x ^k\|.
\]
Rearranging this, we obtain
\[
(\kappa c_k+1)\|u^k - \bar{x}\|^2\leq 2\| c_kg(u^k) - \Pi_K \left( c_kg(u^k) + \lambda^k \right)\|\|\tilde{\lambda}^k - \lambda^k\|+\|u^k-\bar{x}\|^2+\|\bar{x}-x ^k\|^2.
\]
Thus,
\begin{align}\label{6}
\kappa c_k\|u^k - \bar{x}\|^2\leq 2\|c_k g(u^k) - \Pi_K \left(c_k g(u^k) + \lambda^k \right)\|\|\tilde{\lambda}^k - \lambda^k\|+	\|\bar{x}-x ^k\|^2.
\end{align}
From this inequality, we can infer that
\begin{align}\label{claim 1}
\displaystyle\frac{u^k - \bar{x}}{\beta_k} \rightarrow 0,\quad \displaystyle\frac{c_k \|u^k - \bar{x}\|^2}{\beta_k} \rightarrow 0,\quad\mbox{and}\quad \displaystyle\frac{c_k \|u^k - \bar{x}\|^2}{\beta_k^2} \rightarrow 0 \quad \text{as} \quad k \to \infty.
\end{align}
To prove the estimates in (\ref{claim 1}), we use (\ref{6}) together with (\ref{3}),(\ref{5}) and \(c_k \geq \bar{c}\) to obtain
\begin{align*}
\frac{\|u^k - \bar{x}\|^2}{\beta_k^2} &\leq \displaystyle\frac{2}{\kappa c_k}  \frac{\| c_kg(u^k) - \Pi_K \left( c_kg(u^k) + \lambda^k \right)\|}{\beta_k} \cdot \frac{\|\tilde{\lambda}^k - \lambda^k\|}{\beta_k}+ \displaystyle\frac{1}{\kappa c_k} \frac{\|\bar{x}-x^k\|^2}{\beta_k^2}\to 0,\\
\frac{c_k\|u^k - \bar{x}\|^2}{\beta_k} &\leq \displaystyle\frac{2}{\kappa }  \frac{\| c_kg(u^k) - \Pi_K \left( c_kg(u^k) + \lambda^k \right)\|}{\beta_k}  \cdot \|\tilde{\lambda}^k - \lambda^k\|+ \displaystyle\frac{1}{\kappa } \|\bar{x}-x^k\|\cdot\frac{\|\bar{x}-x^k\|}{\beta_k}\to 0,
\end{align*}
and
\begin{align*}
	\frac{c_k\|u^k - \bar{x}\|^2}{\beta_k^2} &\leq \displaystyle\frac{2}{\kappa }  \frac{\| c_kg(u^k) - \Pi_K \left( c_kg(u^k) + \lambda^k \right)\|}{\beta_k} \cdot \frac{\|\tilde{\lambda}^k - \lambda^k\|}{\beta_k}+ \displaystyle\frac{1}{\kappa } \frac{\|\bar{x}-x^k\|^2}{\beta_k^2}\to 0.
\end{align*}

On the other hand, since \( x^k \in B_{\hat{\varepsilon} / (3\hat{l})} (\bar{x}) \) and \( \lambda^k \in B_{\hat{\varepsilon} / (3\hat{l})} (\bar{\lambda}) \) for sufficiently large $k$, and \(u^k\in \mathcal{S}_{c_k}(\lambda^k,x^k)\),  it follows from (\ref{isolated calm for S}) that
\[\|u^k-\bar{x}\|\leq \hat{l}\left(\|x^k-\bar{x}\|+\|\lambda^k-\bar{\lambda}\|\right).
\]
Thus, we conclude that
\begin{align}\label{claim 2}
u^k \to \bar{x}\quad\mbox{as}\quad k \to \infty.
\end{align}

To continue, we apply the first estimate from (\ref{claim 1}) and (\ref{5}) in order to obtain
\[
\xi = \lim_{k \to \infty} \frac{u^k - x^k}{\beta_k} = \lim_{k \to \infty} \frac{u^k - \bar{x}}{\beta_k} - \lim_{k \to \infty} \frac{x^k - \bar{x}}{\beta_k} = 0 - 0 = 0.
\]
Next, we aim to demonstrate that \( \eta = 0 \). When combined with \( \xi = 0 \), this will lead to a contradiction with (\ref{3}). To achieve this, we examine two cases.

First, let's assume that either  \( \{ c_k \} \) or \( \{ c_k / \beta_k \}\) is bounded, then
\begin{align*}
\frac{\|c_{k} g(u^k) - \Pi_K \left(c_{k} g(u^k) +  \lambda^k \right)\|}{\beta_k} &= \frac{c_k}{\beta_k} \|g(u^k) - \Pi_K \left( g(u^k) + c_{k}^{-1} \lambda^k \right) -g(\bar{x})+\Pi_K \left( g(\bar{x}) + c_{k}^{-1} \tilde{\lambda}^k \right)\|\\
&\leq \frac{2\tau_gc_k\| u^k - \bar{x} \|}{\beta_k}   + \frac{\| \tilde{\lambda}^k - \lambda^k \|}{\beta_k},
\end{align*}
where $\tau_g$ represents the Lipschitz constant of $g$ around $\bar{x}$.
It follows from either (\ref{claim 1}) or (\ref{claim 2}) that \( \displaystyle\frac{c_k \| u^k - \bar{x} \|}{\beta_k} \to 0 \). By combining this with (\ref{5}) and taking the limit in the inequality above, we arrive at \( \eta = 0 \), which contradicts (\ref{3}).

Now, suppose that both sequences \( \{ c_k \} \) and \( \{ c_k / \beta_k \} \) are unbounded. Without loss of generality, we can select a subsequence if needed, such that
\[
c_k \to \infty \quad \text{and} \quad \frac{c_k}{\beta_k} \to \infty \quad \text{as} \ k \to \infty.
\]
We denote $z^k:=c_{k} g(u^k)+\lambda^k - \Pi_K \left(c_{k} g(u^k) +  \lambda^k \right)$ for simple. Since \(u^k\in \mathcal{S}_{c^k}(\lambda^k,x^k)\), we deduce from (\ref{optimial_subproblem}) that
\[
0=\nabla_x \mathcal{L}(u^k, \lambda^k, c_k)+c_k^{-1}(u^k-x^k)= \nabla f(u^k)+\nabla g(u^k)z^k+c_k^{-1}(u^k-x^k).
\]
Since \( \tilde{\lambda}^k \in \Lambda(\bar{x})\), this implies that
\[
0=\nabla_x L(\bar{x}, \tilde{\lambda}^k) = \nabla f(\bar{x}) + \nabla g(\bar{x}) \tilde{\lambda}^k.
\]
Combining the last two equations, we have
\begin{align}\label{8}
\nonumber\nabla g(\bar{x})(z^k-\tilde{\lambda}^k)&=
	 \nabla f(\bar{x})- \nabla f(u^k)-\left(\nabla g(u^k)-\nabla g(\bar{x})\right)\left(z^k-\lambda^k\right)-c_k^{-1}(u^k-\bar{x})\\
	 &\qquad\qquad\qquad\qquad\qquad\qquad\qquad+c_k^{-1}(x^k-\bar{x})-\left(\nabla g(u^k)-\nabla g(\bar{x})\right)\lambda^k
\end{align}
On the other hand, it follows from (\ref{2.1}) and $ z^k=c_{k} g(u^k)+\lambda^k - \Pi_K \left(c_{k} g(u^k) +  \lambda^k \right)$  that
\begin{align}\label{9}
	z^k\in \mathcal{N}_{K}\left(g(u^k)-c_k^{-1}(z^k-\lambda^k)\right).
\end{align}
Together with $\tilde{\lambda}^k\in\mathcal{N}_{K}(g(\bar{x}))$ and (\ref{9}), we get from the maximal monotonicity of \(\mathcal{N}_{K}(\cdot) \) in \cite[Corollary 12.18]{RWets98} that
\begin{align}\label{10}
\nonumber0 &\leq \langle z^k - \tilde{\lambda}^k, g(u^k) - c_k^{-1} (z^k - \lambda^k) - g(\bar{x}) \rangle\\
\nonumber&= \langle z^k - \tilde{\lambda}^k, \mathcal{J}g(\bar{x})(u^k - \bar{x}) + O(\| u^k - \bar{x} \|^2) - c_k^{-1}(z^k - \lambda^k)  \rangle\\
&= \langle z^k - \tilde{\lambda}^k, \mathcal{J}g(\bar{x})(u^k - \bar{x}) \rangle - c_k^{-1} \| z^k - \lambda^k \|^2 + \langle z^k - \tilde{\lambda}^k, O(\| u^k - \bar{x} \|^2) \rangle - c_k^{-1} \langle \lambda^k - \tilde{\lambda}^k, z^k - \lambda^k  \rangle.
\end{align}
By combining (\ref{8}) and (\ref{10}), we obtain
\begin{align*}
&\quad\| z^k - \lambda^k \|^2\\ &\leq c_k\langle z^k - \tilde{\lambda}^k, \mathcal{J}g(\bar{x})(u^k - \bar{x}) \rangle  + c_k\langle z^k - \tilde{\lambda}^k, O(\| u^k - \bar{x} \|^2) \rangle - \langle \lambda^k - \tilde{\lambda}^k, z^k - \lambda^k  \rangle\\
&\leq c_k\langle \nabla f(\bar{x})- \nabla f(u^k)-\left(\nabla g(u^k)-\nabla g(\bar{x})\right)\left(z^k-\lambda^k\right)-c_k^{-1}(u^k-\bar{x})+c_k^{-1}(x^k-\bar{x}),u^k - \bar{x} \rangle  \\
&\qquad\qquad-c_k\langle\left(\nabla g(u^k)-\nabla g(\bar{x})\right)\lambda^k,u^k - \bar{x} \rangle + c_k\langle z^k - \tilde{\lambda}^k, O(\| u^k - \bar{x} \|^2) \rangle - \langle \lambda^k - \tilde{\lambda}^k, z^k - \lambda^k  \rangle\\
&\leq c_k\left(O(\| u^k - \bar{x} \|^2)+O(\| u^k - \bar{x} \|^2)\|z^k-\lambda^k\|+O(\| u^k - \bar{x} \|^2)\|\lambda^k\|\right)+\| u^k - \bar{x} \|^2+\|x^k-\bar{x}\|\| u^k - \bar{x} \|\\
&\qquad\qquad\qquad\qquad\qquad\qquad\qquad \qquad\qquad +c_k \|z^k - \tilde{\lambda}^k\|O(\| u^k - \bar{x} \|^2)+\|\lambda^k - \tilde{\lambda}^k\|\|z^k - \lambda^k\|    \\
&\leq (1 + \|z^k-\lambda^k\|+\|z^k - \tilde{\lambda}^k\|+\|\lambda^k\|) O(c_k \| u^k - \bar{x} \|^2) +\| u^k - \bar{x} \|^2+\|x^k-\bar{x}\|\| u^k - \bar{x} \|+ \|\lambda^k - \tilde{\lambda}^k\|\|z^k - \lambda^k\|\\
&\leq (1 + 2\|z^k-\lambda^k\|+\|\lambda^k - \tilde{\lambda}^k\|+\|\lambda^k\|) O(c_k \| u^k - \bar{x} \|^2) +\| u^k - \bar{x} \|^2+\|x^k-\bar{x}\|\| u^k - \bar{x} \|+ \|\lambda^k - \tilde{\lambda}^k\|\|z^k - \lambda^k\|,
\end{align*}
Hence, from the boundedness of $\|\lambda^k\|$, (\ref{3}), (\ref{5}), (\ref{claim 1}) and the fact that
$z^k-\lambda^k =c_{k} g(u^k) - \Pi_K \left(c_{k} g(u^k)+\lambda^k\right)$, we conclude that
\[
\displaystyle\frac{\|c_{k} g(u^k) - \Pi_K \left(c_{k} g(u^k)+\lambda^k\right)\|^2}{\beta_k^2}\rightarrow 0 \quad\mbox{as}\quad k \to \infty.
\]
Apply (\ref{3}) again, we can deduce that
 \( \|\eta\|^2 = 0 \), which implies that
 \( \eta = 0 \). This leads to a contradiction with \( (\xi, \eta) = (0, 0) \), and thus completes the proof. \hfill\textbf{$\Box$}

In the following theorem, we will demonstrate that, under some mild conditions, the inexact PALM given in Algorithm 4.1 automatically satisfies the assumptions (a)-(c) of Theorem 2.2. Consequently, we establish the well-posedness and local convergence of the inexact PALM in Algorithm 4.1.

\begin{theorem}[Well-posedness and Convergence of Inexact PALM]
	Let \((\bar{x}, \bar{\lambda})\) be a KKT point of problem (\ref{1.1}). Suppose that SOSC in (\ref{SOSC}) hold at \((\bar{x}, \bar{\lambda})\) and that \(M_{\bar{x}}\) is calm at the origin for \(\bar{\lambda}\). Then, positive constants \(\hat{\alpha}\), \(\bar{c}\), and \(\varepsilon_0\) can be found such that for any initial point \((x^0, \lambda^0) \in B_{\varepsilon_0}(\bar{x}, \bar{\lambda})\) and any sequence \(\{c_k\}\) with \(c_k \geq \bar{c}\), there exists a sequence \(\{(x^k, \lambda^k)\}\) satisfying the update rules (\ref{x-pupdate}) and (\ref{lambda-update}) with \(\epsilon_k = o(r(x^k, \lambda^k))\), and fulfilling (\ref{esta}), where \(r\) denotes the residual function defined in (\ref{residual}).
	Furthermore, any such sequence \(\{(x^k, \lambda^k)\}\) converges to \((\bar{x}, \hat{\lambda})\) for some \(\hat{\lambda} \in \Lambda(\bar{x})\). Both the sequence \(\{(x^k, \lambda^k)\}\) converges to \((\bar{x}, \hat{\lambda})\) and the distance \(\{\text{dist}((x^k, \lambda^k), \{\bar{x}\} \times \Lambda(\bar{x}))\}\) converges to zero at a Q-linear rate. Additionally, if \(c_k \to \infty\), the convergence rates for both sequences become Q-superlinear.
\end{theorem}
\textbf{Proof.} To prove the results, we check the validity of assumptions (a)--(c) of Theorem 2.2 as follows. We begin with assumption (a). Since $\mathcal{S}_{KKT}$ represents the solution mapping for the canonical perturbation of (\ref{GE3}), it suffices to confirm the calmness of the solution mapping \(\mathcal{S}_{KKT}\)  at the origin for \((\bar{x}, \bar{\lambda})\). It follows from Proposition 2.2 that the calmness of $M_{\bar{x}}$ at origin for $\bar{\lambda}$, together with the SOSC in (\ref{SOSC}), ensures the semi-isolated calmness of the solution map $\mathcal{S}_{KKT}$ at origin for $(\bar{x},\bar{\lambda})$, hence  \(\mathcal{S}_{KKT}\) is calm at the origin for \((\bar{x}, \bar{\lambda})\).

Next, we proceed to assumption (b) in Theorem 2.2. The concept of the semi-isolated calmness property of  \(\mathcal{S}_{KKT}\) in Definition 2.3 also implies  the existence of  a  constant $\delta_1 >0$ such that
\begin{align*}
	\mathcal{S}_{KKT}(0,0)\cap B_{\delta_1}(\bar{x}, \bar{\lambda})\subset\left( \{\bar{x}\} \times  \Lambda(\bar{x})\right)\cap B_{\delta_1}(\bar{x}, \bar{\lambda} ),
\end{align*}
which means that
\begin{align*}
\|\tilde{x} - \bar{x}\| + \text{dist}(\tilde{\lambda}, \Lambda(\bar{x})) \leq \text{dist}\left((\tilde{x}, \tilde{\lambda}),\mathcal{S}_{KKT}(0,0)\right)
\end{align*}
holds for any $(\tilde{x}, \tilde{\lambda})\in B_{\delta_1}(\bar{x}, \bar{\lambda})$. Combined  with (\ref{r}), there exist positive constants $\delta_2\leq \delta_1$ and $\hat{\kappa}$ such that for any $(\tilde{x}, \tilde{\lambda})\in B_{\delta_2}(\bar{x}, \bar{\lambda})$,
\begin{align}\label{11}
r(\tilde{x}, \tilde{\lambda})\leq \hat{\kappa}\left(\|\tilde{x} - \bar{x}\| + \text{dist}(\tilde{\lambda}, \Lambda(\bar{x})) \right)\leq \hat{\kappa}\text{dist}\left((\tilde{x}, \tilde{\lambda}),\mathcal{S}_{KKT}(0,0)\right).
\end{align}
On the other hand, by Proposition 3.2, we can find constants $\delta_3>0$ and $\hat{\alpha}>0$ such that for any $c > \bar{c}$ and $(\tilde{x}, \tilde{\lambda}) \in B_{\delta_3}(\bar{x}, \bar{\lambda})$ satisfying $r(\tilde{x}, \tilde{\lambda}) > 0$ and the optimal solution \( x\in \mathcal{S}_{c}(\tilde{\lambda},\tilde{x}) \), the following inequality holds:
					\begin{equation}\label{15}
						\|x - \tilde{x}\| + \|\lambda - \tilde{\lambda}\| \leq \hat{\alpha} r(\tilde{x}, \tilde{\lambda}),
					\end{equation}
where  $\lambda:=c g(x)+\tilde{\lambda} - \Pi_K \left( c g(x) +  \tilde{\lambda} \right)$. Take $\varepsilon_1:=\min{(\delta_2,\delta_3)}$. Employing (\ref{11}), for any  $(\tilde{x}, \tilde{\lambda}) \in B_{\varepsilon_1}(\bar{x}, \bar{\lambda})$, it holds that
\begin{align}\label{12}
	\|x - \tilde{x}\| + \|\lambda - \tilde{\lambda}\|  \leq \hat{\alpha}\hat{\kappa}\text{dist}\left((\tilde{x}, \tilde{\lambda}),\mathcal{S}_{KKT}(0,0)\right) .
\end{align}
If \((\tilde{x}, \tilde{\lambda}) \in B_{\varepsilon_1}(\bar{x}, \bar{\lambda})\) and \(r(\tilde{x}, \tilde{\lambda}) = 0\), then $(\tilde{x}, \tilde{\lambda})\in \mathcal{S}_{KKT}(0,0)$. Let $x:=\tilde{x}$ and $\lambda:=\tilde{\lambda}$, in which case (\ref{12}) then holds trivially.
Therefore, we conclude that for any $(\tilde{x}, \tilde{\lambda}) \in B_{\varepsilon_1}(\bar{x}, \bar{\lambda})$, there exists a pair \((x, \lambda)\in X\times H\) such that
\[
0 \in \mathcal{A}(x, \lambda,\tilde{x}, \tilde{\lambda}, \epsilon, c) + G(x, \lambda),
\]
and the estimate (\ref{12}) holds,
where \(\mathcal{A}(\cdot)\) and \(G(\cdot)\) are defined in (\ref{A1}) and (\ref{GE2}) respectively.

Finally, we justify assumption (c) in Theorem 2.2.  Consider any point \((\tilde{x}, \tilde{\lambda}) \in B_{\varepsilon_1}(\bar{x}, \bar{\lambda})\), and any point \((x, \lambda) \in X \times H\) such that \(\| (x, \lambda) - (\tilde{x}, \tilde{\lambda}) \| \leq \hat{\kappa}\hat{\alpha} \mbox{dist}((\tilde{x}, \tilde{\lambda}) , \mathcal{S}_{KKT}(0, 0))\). Additionally, let \(c \geq \bar{c}\), and choose  \(\epsilon\) such that \(
\epsilon = o(r(\tilde{x}, \tilde{\lambda}))\). For any \((w_1, w_2) \in \Phi(x, \lambda) - \mathcal{A}(x, \lambda, \tilde{x}, \tilde{\lambda}, \epsilon, c)\), where \(\Phi\) is defined in (\ref{GE2}), and combining with (\ref{11}) and (\ref{15}), we can find an element $e\in X$ with $\|e\|=1 $ such that
\begin{align*}
\|(w_1, w_2)\| &= \|\left(\epsilon e + c^{-1}(x-\tilde{x}), c^{-1}(\lambda- \tilde{\lambda})\right)\|\\
&\leq \epsilon + c^{-1}\left(\|x-\tilde{x}\|+  \| \lambda - \tilde{\lambda} \|\right) \\
&\leq o(r(\tilde{x}, \tilde{\lambda})) + c^{-1}\hat{\alpha}r(\tilde{x}, \tilde{\lambda})\\
&=	\left(\frac{o(r(\tilde{x}, \tilde{\lambda}))}{r(\tilde{x}, \tilde{\lambda})}+c^{-1}\hat{\alpha}\right)\cdot r(\tilde{x}, \tilde{\lambda}) \\
&\leq \left(\frac{\hat{\kappa}o(r(\tilde{x}, \tilde{\lambda}))}{r(\tilde{x}, \tilde{\lambda})}+c^{-1}\hat{\alpha}\hat{\kappa}\right)\mbox{dist}((\tilde{x}, \tilde{\lambda}) , \mathcal{S}_{KKT}(0, 0)).
\end{align*}
We define \(\omega : X \times H\ \times X \times H  \times \mathbb{R}_+ \times (0, +\infty) \to \mathbb{R}_+\) by
\[
\omega(x, \lambda, \tilde{x}, \tilde{\lambda}, \epsilon, c) = \frac{\hat{\kappa} o(r(\tilde{x}, \tilde{\lambda}))}{r(\tilde{x}, \tilde{\lambda})}+c^{-1}\hat{\kappa}\hat{\alpha}.
\]
Thus, we have established the second inequality in assumption (c) of Theorem 2.2. Since \(\omega(x, \lambda, \tilde{x}, \tilde{\lambda}, \epsilon, c) \to 0\) as \((x, \lambda, \tilde{x}, \tilde{\lambda},\epsilon, c) \to (\bar{x}, \bar{\lambda}, \bar{x}, \bar{\lambda}, 0,\infty)\), the first part in assumption (c) is also satisfied .
				
Finally, applying Theorem 2.2, we establish all the convergence results, thereby completing the proof.
 \hfill\textbf{$\Box$}
 \begin{remark}
  When \( X := \mathbb{R}^n \) and \( K := \mathbb{R}_-^m \), the problem (\ref{1.1}) reduces to a NLP problem. In this case, the calmness of \( M_{\bar{x}} \) at the origin for \( \bar{\lambda} \) is always guaranteed and can be established through the classical Hoffman’s lemma (see, e.g., \cite{DR2014}). Therefore, the conclusions in Theorem 3.1 for the NLP case hold only under the second-order sufficient conditions, which are consistent with the results in \cite[Theorem 4.3]{IS2024}. Thus, Theorem 3.1 extends the convergence results of the PLAM for the NLP case in \cite{IS2024} to a broader class of conic optimization problems with non-polyhedral constraints.

  \end{remark}
\section{Conclusion}
\setcounter{equation}{0}
In this study, we analyze the local convergence properties and convergence rates of the Proximal Augmented Lagrangian Method (PALM) for a class of nonconvex conic optimization problems under weaker assumptions than those in existing literature. We rigorously establish that both the exact and inexact versions of PALM exhibit Q-linear convergence rates under second-order sufficient optimality conditions and the calmness of the multiplier function. Furthermore, we demonstrate that the convergence rate becomes superlinear when the penalty parameter $c_k$ increases to $+\infty$. As a byproduct, we provide a comprehensive proof of the solvability for PALM subproblems and establish the uniform isolated calmness property of the optimal solution mapping for the subproblem.
 
\end{document}